\documentclass[conference]{IEEEtran}
\IEEEoverridecommandlockouts
\usepackage{amsmath,amssymb,amsfonts}
\usepackage{algorithmic}
\usepackage{graphicx}
\usepackage{textcomp}
\usepackage{xcolor}
\usepackage{hyperref}
\usepackage{cite}
\usepackage{bbm}
\usepackage{subcaption}

\usepackage[ruled]{algorithm2e}

\begin{document}

\title{Learned Finite-Time Consensus for Distributed Optimization}

\author{\IEEEauthorblockN{Aaron Fainman and Stefan Vlaski}
\thanks{The authors are with the Department of Electrical and Electronic Engineering, Imperial College London. Emails: \{aaron.fainman22, s.vlaski\}@imperial.ac.uk.}
}

\maketitle

\begin{abstract}
Most algorithms for decentralized learning employ a consensus or diffusion mechanism to drive agents to a common solution of a global optimization problem. Generally this takes the form of linear averaging, at a rate of contraction determined by the mixing rate of the underlying network topology. For very sparse graphs this can yield a bottleneck, slowing down the convergence of the learning algorithm. We show that a sequence of matrices achieving finite-time consensus can be learned for unknown graph topologies in a decentralized manner by solving a constrained matrix factorization problem. We demonstrate numerically the benefit of the resulting scheme in both structured and unstructured graphs.
\end{abstract}


\section{Introduction and Related Work}
We consider the problem of distributed optimization, in which a network of $K$ agents collectively aim to solve an aggregate optimization problem of the form:
\begin{equation}\label{eq:consensus_opt}
    \min_{w\in\mathbb{R}^M} \quad J(w)=\frac{1}{K}\sum_{k=1}^K J_k(w)
\end{equation}
Here, each local objective function $J_k: \mathbb{R}^M\rightarrow\mathbb{R}$ is known only to agent $k$. The problem arises frequently in settings such as smart grids, wireless sensor networks, and large scale machine learning (see, \textit{e.g.}~\cite{mateos,kekatos, assran, vlaski, sayed_book }).

Decentralized methods of the diffusion~\cite{cattivelli, sayed_book} or consensus~\cite{ram} type converge to the optimal solution under diminishing step-sizes, at the cost of convergence rate. When employing a constant step-size and exact gradients, they exhibit a non-negligible bias. This bias can be corrected with algorithms such as ADMM~\cite{admm}, EXTRA~\cite{extra}, Exact Diffusion \cite{exact_diff}, gradient tracking (GT)~\cite{next} or AugDGM~\cite{augdgm}. The latter two correct the steady-state bias using an auxiliary variable that tracks the average gradient. Many variants of each of these algorithms exist, based on whether the graph is directed, whether the network is time-varying, and whether the communication is asynchronous~\cite{diging, add_opt, asynch}.
As a representative example, the AugDGM recursion~\cite{augdgm} is shown in~\eqref{eq:aug1_node} and~\eqref{eq:aug2_node}, using the local optimization variable \( w_{k}^{(i)} \), auxiliary gradient tracking variable $g_{k, i}$, step-size $\mu > 0$, and combination weights at iteration $i$, $a_{k \ell}^{(i)} = [A^{(i)}]_{k, \ell}$:
\begin{align}
    w^{(i)}_k &= \sum_{\ell\in\mathcal{N}_k} a^{(i)}_{k\ell}\left(w^{(i-1)}_{\ell}-\mu g_\ell^{(i-1)}\right)  \label{eq:aug1_node}\\
    g^{(i)}_k &= \sum_{\ell\in\mathcal{N}_k} a^{(i)}_{k\ell}\left(g^{(i-1)}_{\ell} + \nabla J_{\ell}(w^{(i)}_{\ell})- \nabla J_{\ell}(w^{(i-1)}_{\ell}) \right)  \label{eq:aug2_node}
\end{align}
Here, \( \mathcal{N}_k \) denotes the neighborhood of agent \( k \), which is the set of agents with whom agent \( k \) is able to communicate. We assume that the combination weights are consistent with the underlying network topology, meaning that \( a^{(i)}_{k\ell} = 0 \) if \( \ell \notin \mathcal{N}_k \). Relation~\eqref{eq:aug1_node} corresponds to an ATC-diffusion~\cite{sayed_book} type recursion along the estimate of the aggregate gradient \( g_\ell^{(i-1)} \), while~\eqref{eq:aug2_node} updates the \( g_{k}^{(i)} \) following a dynamic consensus mechanism~\cite{next}. We can write the recursions more compactly by defining network level quantities $w^{(i)}\triangleq \mathrm{col}\{w^{(i)}_{1},\dots,w^{(i)}_{K}\}$, with the gradient of $J(w)$ given by $\nabla J(w^{(i)}) = \mathrm{col}\{\nabla J_1(w_1^{(i)}),\dots,\nabla J_{K}(w_K^{(i)})\}$ and the combination weights, $\mathcal{A}^{(i)}\triangleq A^{(i)}\otimes I_M$~\cite{augdgm}:
\begin{align}
    w^{(i)} &= \mathcal{A}^{(i)}\left(w^{(i-1)}-\mu g^{(i-1)}\right) \label{eq:aug1} \\
    g^{(i)} &= \mathcal{A}^{(i)}\left(g^{(i-1)}+\nabla J(w^{(i)})-\nabla J(w^{(i-1)})\right) \label{eq:aug2}
\end{align}
When tackling aggregate optimization problems of the form~\eqref{eq:consensus_opt}, the performance of most decentralized learning algorithms is limited by two types of errors, the optimization error and the consensus error, which is defined as the variation between individual local variables.
Some algorithms perform multiple rounds of communication between each local computation to mitigate this variation (see \textit{e.g.}~\cite{jakovetic}). This comes at the cost of increased communication complexity. Two recent algorithms have improved the convergence rate by reducing the consensus error using a loopless Chebyshev acceleration technique and by incorporating momentum terms~\cite{optimalgt, huang}. 

The consensus error can also be reduced by optimizing the network's combination weights, collected in the matrix \( [A^{(i)}]_{k, \ell} \triangleq a_{k \ell}^{(i)} \). This has been investigated in the context of linear averaging. To motivate this procedure, note that most decentralized optimization algorithms, in one form or another, contain a mixing step of the form:
\begin{align}
    x^{(i)} = \mathcal{A}^{(i)} x^{(i-1)} \label{eq:lin-avg} 
\end{align}
where $\mathcal{A}^{(i)}$ is the combination matrix at iteration $i$, and $x^{(i)}\in\mathbb{R}^{KM\times 1}$ contains local estimates of some quantity to be averaged. It is well established that the iterates will asymptotically converge to the average of the local quantities at initialization, provided that each combination matrix in the sequence is doubly-stochastic and primitive, with spectral radius equal to one~\cite{xiao, hendrickx}. 

In~\cite{xiao}, the convergence of linear averaging is improved by optimizing a static combination matrix \( \mathcal{A}^{(i)} \triangleq \mathcal{A} \) to have the fastest mixing rate. For symmetric matrices, this is equivalent to minimizing the spectral norm, $\lVert A- \frac{1}{K}\mathbbm{1}\mathbbm{1}^\mathsf{T}\rVert_2$, which can be solved with a semi-definite program, subject to appropriate constraints on $A$. Results show improved convergence over traditional weighting rules \cite{xiao}.

If time-varying combination matrices are allowed, the exact average can be found in finite time. Consider a sequence of doubly-stochastic weight matrices. If the sequence forms a sparse factorization of the scaled all-ones matrix, such that:
\begin{equation}
    A_{\tau}\cdots A_2A_1=\frac{1}{K}\mathbbm{1}\mathbbm{1}^\mathsf{T}
\end{equation}
then the exact average can be found in $\tau$ steps \cite{Cesar, kibangou, safavi, sandryhaila}. The quantity $\tau$ is known as the graph's consensus number and is bounded below by the graph diameter and above by twice the graph radius \cite{hendrickx}. For certain families of graphs, exact sequences of finite-time (FT) matrices are known. This includes hypercubes, 1-peer exponential graphs \cite{Cesar}, and p-peer hyper-cuboids \cite{Cesar}. In \cite{Cesar}, FT matrices are investigated for both linear averaging and more general optimization problems of the form~\eqref{eq:consensus_opt}. The resulting GT algorithm with FT sequences achieves identical iteration complexity as a static weighting rule, despite having much sparser communication links at each iteration, resulting in improved communication cost.

The discussion thus far motivates the use of time-varying combination matrices to improve convergence of decentralized optimization algorithms. These combination weights require specific network topologies that are known beforehand. Other works like, \cite{xiao} and in \cite{sparsefactors}, are centralized in their implementation. Motivated by these limitations, we develop a \emph{decentralized} algorithm for learning optimal combination sequences for linear averaging, which do not require prior knowledge of the network topology. Numerical demonstrate that the learned combination sequences lead to improved decentralized optimization even when they are found only approximately.

\section{Algorithm Development}
\noindent We aim to find the sequence of symmetric, doubly-stochastic weight matrices, $\{A_j\}_{j=1}^{\tau}$ over a graph $\mathcal{G}=(\mathcal{V}, \mathcal{E})$ whose product equals the scaled all-ones matrix:
\begin{align}
    \min_{A_{\tau}, ..., A_2A_1}\quad & \left|\left| \frac{1}{K}\mathbbm{1}\mathbbm{1}^\mathsf{T} - A_{\tau}\cdots A_2A_1\right|\right|_F^2 \label{eq:main_prob}\\ 
    \text{subject to}\quad & A_j = A_j^\mathsf{T} \label{eq:symmetric}\\
    & [A_j]_{k,\ell} = 0 \ \quad\ \text{if}\quad (k,\ell)\notin\mathcal{E} \label{eq:sparse}\\
    & A_j \mathbbm{1} =\mathbbm{1} \quad \text{for}\quad j\in\{1,...,\tau\} \label{eq:stochastic}
\end{align}
where $||\cdot||_F$ is the Frobenius norm. Owing to the non-convex objective function, this is a non-trivial problem for $\tau>1$ \cite{lnn}. 

We pursue the minimization of \eqref{eq:main_prob} using a projected gradient descent algorithm. At each iteration, for each matrix $\{A_j\}_{j=1}^\tau$, we take a gradient step to obtain an intermediate estimate $A'^{(i)}_j$:
\begin{align}\label{eq:pgd:grad}
    A^{'(i)}_j =A^{(i-1)}_j-\nabla_{A_j} J(A^{(i-1)}_1, \ldots, A^{(i-1)}_{\tau})
\end{align}
Then, we project back onto the feasible set defined by the constraints in~\eqref{eq:main_prob}--\eqref{eq:stochastic}:
\begin{align}\label{eq:pgd:proj}
    A^{(i)}_j = P_{ \mathcal{C}_{\mathrm{sym}} \cap \mathcal{C}_{\mathrm{spars}} \cap \mathcal{C}_{\mathrm{sym}}}\left(A^{'(i)}_j\right)
\end{align}
Here, \( \mathcal{C}_{\mathrm{sym}} \) denotes the set of symmetric matrices~\eqref{eq:symmetric}, \( \mathcal{C}_{\mathrm{spars}} \) denotes the set of matrices whose sparsity pattern aligns with the network's underlying topology~\eqref{eq:sparse} and \( \mathcal{C}_{\mathrm{stoc}} \) denotes the set of right stochastic matrices~\eqref{eq:stochastic}. Together, these conditions ensure that employing the learned sequence of matrices in a decentralized learning algorithm such as~\eqref{eq:aug1}--\eqref{eq:aug2} will result in convergence to the minimizer of~\eqref{eq:consensus_opt}, while being amenable to a decentralized implementation. The gradient step in \eqref{eq:pgd:grad} for matrix $j$ at iteration $i$ is given by:
\begin{equation}\label{eq:pgd:actual_grad}
    A'^{(i)}_j = A^{(i-1)}_j-\mu A^{(i-1)}_{j+1:\tau} A^{(i-1)}_{\tau:j+1} A^{(i-1)}_j A^{(i-1)}_{j-1:0}A^{(i-1)}_{0:j-1}
\end{equation}
where we define the product, $A_{j:k}\triangleq A_j\dots A_{k}$ for brevity.

A decentralized closed-form expression for the projection onto the intersection \( \mathcal{C}_{\mathrm{sym}} \cap \mathcal{C}_{\mathrm{spars}} \cap \mathcal{C}_{\mathrm{sym}} \) does not exist. Instead, we replace the projection step in~\eqref{eq:pgd:proj} with sequential projections so that:
\begin{align}\label{eq:pgd:seq_proj}
    A^{(i)} =  P_{\mathcal{C}_{\mathrm{sym}}} \left( P_{\mathcal{C}_{\mathrm{stoc}} \cap \mathcal{C}_{\mathrm{spars}}}\left(P_{\mathcal{C}_{\mathrm{spars}}}(A^{'(i)}) \right) \right)
\end{align}
The projection onto the set of sparse matrices is straightforward, as it corresponds to nulling the entries \( [A]_{k, \ell} \) corresponding to missing edges \( (k, \ell) \notin \mathcal{E} \). The joint projection onto the set of right stochastic and sparse matrices corresponds to: 
\begin{align}
\left[P_{\mathcal{C}_{\mathrm{stoc}} \cap \mathcal{C}_{\mathrm{spars}} }(A)\right]_{k, \ell} = 
\begin{cases}\label{eq:proj:stoch}
    [A]_{k,\ell} - \frac{\sum_{j\in\mathcal{N}_k}([A]_{k,j})-1}{\lvert \mathcal{N}_k\rvert}, & \textrm{if} (k,\ell)\in\mathcal{E}, \\
      0, & \text{otherwise.}
\end{cases}
\end{align}
\normalsize
Here, \( |\mathcal{N}_k| \) is the size of the neighborhood of agent \( k \). This is also a local operation, given that \( \sum_{j\in\mathcal{N}_k}[A]_{k,j} \) is the sum of the non-zero elements of row $k$ in $A$, or the sum of the agent's weights. The symmetry projection is given by:
\begin{align}\label{eq:proj:sym}
    [P_{\mathcal{C}_{\mathrm{sym}}}(A)]_{k, \ell} = \frac{1}{2}([A]_{k,\ell} + [A]_{\ell,k})
\end{align}
which requires one round of communication between an agent and its neighbors. 

The fact that the projections in \eqref{eq:pgd:seq_proj} are computed sequentially, rather than projecting onto the intersection of all constraints, means that $A_j^{(i)}$ may not be in the feasible set of \eqref{eq:main_prob}--\eqref{eq:stochastic} during the transient stage of the algorithm. Even in the convex setting, feasibility of the iterates of the incrementally projected gradient algorithm is only guaranteed asymptotically~\cite{pierro}. For example, performing the projection onto the set of right-stochastic matrices~\eqref{eq:proj:stoch}, and then the projection onto the set of symmetric matrices~\eqref{eq:proj:sym}, results in a matrix that is symmetric but not necessarily right stochastic. 

In order to ensure that, even when stopping after a finite number of iterations, the resulting sequence of matrices is suitable for decentralized optimization, we need to make a suitable correction of the final iterate. The order of projections in~\eqref{eq:pgd:seq_proj} ensures that \( A_j^{(i)} \) is symmetric for all \( i \), since the symmetry projection is performed last. It will also be sparse, since the symmetry projection~\eqref{eq:proj:sym} preserves the sparsity pattern of the combination policy. We are however not guaranteed that at any finite time instant \( i \), \( A_j^{(i)} \) is right-stochastic. We can ensure this post-hoc. At the final iteration, the diagonal entries of \( A_j^{(T)} \), are adjusted as follows:
\begin{align}\label{eq:diag}
    [A_j^{(T)}]_{k,k} = 1-\sum_{\ell\in\mathcal{N}_k}[A_j^{(T)}]_{k,\ell}
\end{align}
 thereby ensuring that \( A_j^{(T)} \) is in the feasible set defined by~\eqref{eq:symmetric}--\eqref{eq:stochastic}. It is known that projected gradient descent with sequential projections converges for convex problems \cite{pierro}. Owing to the non-convex structure of~\eqref{eq:main_prob}, proving convergence is not guaranteed by this result.
 Indeed, we can interpret~\eqref{eq:main_prob} as the optimization problem for a linear deep neural network, subject so some structural constraints which make the resulting weights appropriate for decentralized optimization. Despite their non-convex nature, the weights of linear neural networks have been shown to converge linearly under suitable conditions on the initialization and the target function~\cite{lnn}. While our setting does not fully match the conditions of~\cite{lnn}, our numerical results nevertheless indicate that Algorithm 1 does converge to an optimal solution. Establishing conditions for linear convergence is left for future work.

The full Learned Finite-Time Consensus (LFTC) algorithm is shown in Algorithm~\ref{alg:pgd}.  
The algorithm requires agents to agree a priori on the sequence length \( \tau \). Ideally, \( \tau \) should be set to the consensus number of the graph, which is of course generally unknown to the network, and especially to individual agents. As we demonstrate in our numerical results, however, the proposed scheme is robust to the choice of \( \tau \), and even a misspecified consensus number yields substantial improvements in the rate of convergence.

\begin{algorithm}
    \SetAlgoLined
    \caption{Learned Finite-Time Consensus (LFTC)}\label{alg:pgd}
        Denote $A^{(i)}_{p:q} = A^{(i)}_p\dots A^{(i)}_{q}$ \\
        \textbf{Initialize} $A^{(0)}_j\in \mathbb{R}^{K\times K}\quad \forall\  j\in\{1\dots\tau\}$ \\
        \For{$T$ \text{iterations}}{
        \For{$j\gets 1$ \KwTo $\tau$}{
        ${ A^{'(i)}_j =  A^{(i-1)}_j-\mu A^{(i-1)}_{j+1:\tau} A^{(i-1)}_{\tau:j+1} A^{(i-1)}_j A^{(i-1)}_{j-1:1}A^{(i-1)}_{1:j-1} }$\\
        $A^{(i)} =  P_{\mathcal{C}_{\mathrm{sym}}} \left( P_{\mathcal{C}_{\mathrm{stoc}} \cap \mathcal{C}_{\mathrm{spars}}}\left(P_{\mathcal{C}_{\mathrm{spars}}}(A^{'(i)}) \right) \right)$\\
        }}
        \For{$j\gets 1$ \KwTo $\tau$}{
        \For{$k\gets 1$ \KwTo $K$}{
        $[A^{(T)}_j]_{k,k} = 1-\sum_{\ell\in\mathcal{N}_k}[A^{(T)}_j]_{k,\ell}$
        }
        }
    \end{algorithm}
Having learned the sequence of \( \tau \) combination matrices in the LFTC algorithm, we can then incorporate these into an averaging scheme of the form~\eqref{eq:lin-avg} by setting:
\begin{align}\label{eq:lftc_averaging}
    A^{(i)} \triangleq A_{i \% \tau + 1}^{\mathrm{LFTC}}
\end{align}
where \(A_{j}^{\mathrm{LFTC}}\) refers to the $j^{th}$ matrix in the sequence returned by the LFTC algorithm and \( \% \) denotes the modulo operation. If \( \tau \) corresponds to the consensus number of the graph, and the collection \( \{ A_j^{(T)} \}_{j=1}^{\tau} \) minimizes~\eqref{eq:main_prob} exactly, then~\eqref{eq:lin-avg} with the choice~\eqref{eq:lftc_averaging} will reach exact consensus after exactly \( \tau \) iterations. Otherwise, the algorithm will cycle through the sequence of combination matrices until convergence. Similarly, we can use the choice~\eqref{eq:lftc_averaging} in recursions~\eqref{eq:aug1}--\eqref{eq:aug2} to yield an LFTC-variance of AugDGM. We evaluate the performance of both of these schemes numerically in the next section.

\section{Experimental Results}
\begin{center}
\begin{figure}[]
    \centering 
\begin{subfigure}{0.85\linewidth}
  \includegraphics[width=\linewidth]{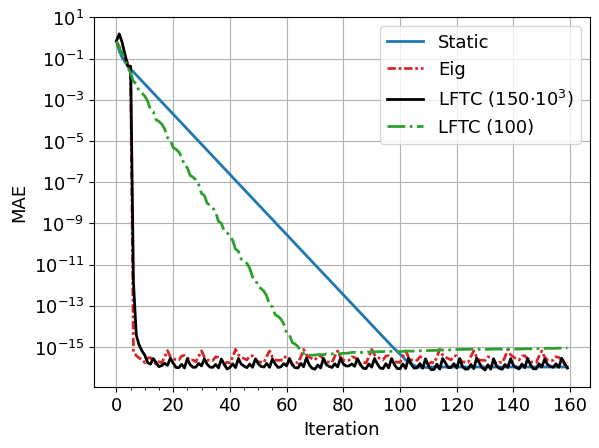}
\end{subfigure}
\caption{Linear averaging over a hypercube with 64 agents. }
\label{fig:results_hyp}
\end{figure}
\end{center}

We demonstrate the benefit of the sequences~\eqref{eq:lftc_averaging} obtained via the LFTC Algorithm~\ref{alg:pgd} for two examples of the averaging recursion~\eqref{eq:lin-avg}, as well as in the AugDGM recursions~\eqref{eq:aug1}--\eqref{eq:aug2}. The LFTC-based algorithms are compared against a static combination matrix, given by the fastest mixing rate in \cite{xiao}, as well as exact finite-time consensus sequences generated with full knowledge of the graph. There are a number of techniques for generating an exact finite-time consensus sequence for a given graph~\cite{sandryhaila, safavi, kibangou}. We use as a baseline the construction ${A_s = \frac{1}{1-\lambda_s}(W-\lambda_sI_K)}$ where $\lambda_s$ are the distinct eigenvalues of $W$ less than one, and $W$ is a doubly-stochastic combination matrix consistent with the graph topology. It bears emphasising that both the fastest mixing static combination rule, and the FT consensus sequence require full central knowledge of the graph topology, and are merely used as benchmarks. Machine precision for the results in Fig.~\ref{fig:results_hyp} is $10^{-15}$.

Fig.~\ref{fig:results_hyp} shows the mean absolute error (MAE) for the linear averaging recursion~\eqref{eq:lin-avg} for a hypercube with 64 nodes. The proposed scheme is shown twice: after running the LFTC algorithm until convergence ($150\cdot10^3$ iterations), and following early stopping after $100$ iterations to illustrate whether improvements are still observed when~\eqref{eq:main_prob} is solved only approximately. The eigenvalue rule and the full LFTC algorithm have identical performance for the linear averaging, and converge after only $6$ iterations, which corresponds to the consensus number of the hypercube. Clearly, the LFTC descent algorithm is able to learn the FT consensus sequence without any prior knowledge of the structure of the graph. The performance of the $100$-iteration LFTC algorithm does degrade, but still shows significant improvement over the static fastest mixing rate combination rule. 
\begin{center}
\begin{figure}[htb]
    \centering 
\begin{subfigure}{0.85\linewidth}
  \includegraphics[width=\linewidth]{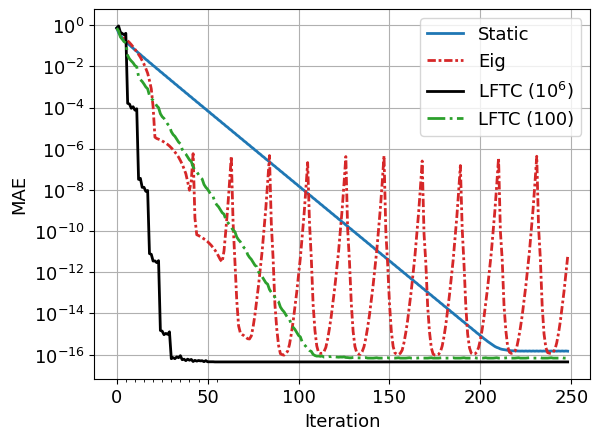}
\end{subfigure}
\caption{ Linear averaging over a Barabasi-Albert model with 32 agents.}
\label{fig:results_pag}
\end{figure}
\end{center}

Results when averaging over a 32-node Barabasi-Albert model~\cite{pag} are shown in Fig.~\ref{fig:results_pag}. Compared to the hypercube, a Barabasi-Albert graph has significantly less symmetry, and a stronger community structure. The graph has a radius of 3 and so we choose a sequence length \( \tau \) of 6 matrices for the LFTC. We note that compared to the hypercube, the LFTC algorithm converges more slowly to the optimal learned weights for a Barabasi-Albert graph, resulting in an error of \( 0.1\% \) between the matrix product $A_\tau \dots A_2 A_1$ and the scaled all-ones matrix after $10^{6}$ iterations. This prevents the resulting averaging recursion~\eqref{eq:lin-avg} from converging after \( 6 \) iterations. Nevertheless in the linear averaging recursion, we observe significant speedup compared to a static combination rule. The $100$-iteration LFTC is slower but still outperforms static averaging. The eigenvalue rule results in significantly less stable behavior. We attribute this to a combination of numerical errors in the eigendecomposition, resulting in inexact computation of the consensus sequence, coupled with the fact that each individual \( A_s \) in the sequence is not necessarily constrained to being stable.
\begin{center}
\begin{figure}[htb]
    \centering 
\begin{subfigure}{0.85\linewidth}
  \includegraphics[width=\linewidth]{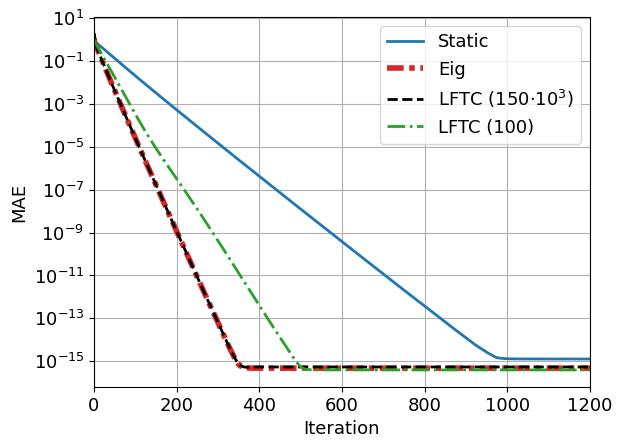}
\end{subfigure}
\caption{ MAE for a linear least squares problem on a hypercube. }
\label{fig:results_aug}
\end{figure}
\end{center}

Performance for a linear regression problem with a least squares objective \( J_k(w) = \frac{1}{2N} \sum_{n=1}^N (\gamma_{k, n} - h_{k, n}^{\mathsf{T}} w)^2 \) on a 64-node hypercube is shown in Fig.~\ref{fig:results_aug}. Data is generated using a linear model \( \boldsymbol{\gamma}_{k} = \boldsymbol{h}_{k}^{\mathsf{T}} w^o + \boldsymbol{v}_{k}\) where the features \( \boldsymbol{h}_k \) follow a normal distribution $\mathcal{N}(0,I_M)$ and \( M = 20 \). The noise is also normally distributed with \( \sigma_v^2 = 0.1\). Optimization follows the AugDGM algorithm in~\eqref{eq:aug1} and~\eqref{eq:aug2}. Once again, the eigenvalue rule and fully-converged LFTC algorithm have identical performance, converging after 350 iterations. The $100$-iterate LFTC algorithm converges in 500 iterations, which is still significantly faster than the static combination policy.

The simulations indicate that the LFTC algorithm is able to learn effectivly finite-time consensus sequences for a graph. Results also demonstrate the robustness of the algorithm. Even when the the sequence does not converge exactly to $\frac{1}{K}\mathbbm{1}\mathbbm{1}^\mathsf{T}$, there is a benefit to the performance for decentralized optimization. 

\section{Conclusion}
We have developed the Learned Finite-Time Consensus (LFTC) Algorithm for finding sequences of combination matrices for decentralized optimization. In contrast to prior works, the scheme does not require prior knowledge of the graph topology, and instead learns an optimal sequence of mixing matrices by solving a sparse matrix factorization problem. Numerical results demonstrate that the resulting sequences, even when found only approximately, result in substantial acceleration of decentralized optimization algorithms.

\bibliographystyle{IEEEtran}
\bibliography{References}

\end{document}